\documentclass[12pt]{amsart}

\usepackage{amssymb,amsthm,amsfonts,graphicx,comment}
\usepackage[usenames]{color}
\usepackage{verbatim}

\newtheorem{thm}{Theorem}[section]
\newtheorem{lemma}[thm]{Lemma}
\newtheorem{cor}[thm]{Corollary}
\newtheorem{prop}[thm]{Proposition}

\newtheorem*{sublemma*}{Sublemma}

\theoremstyle{definition}
\newtheorem{definition}[thm]{Definition}
\newtheorem{remark}[thm]{Remark}
\newtheorem{warning}[thm]{Warning}

\newtheorem{example}[thm]{Example}

\DeclareMathOperator{\genus}{\mathsf{genus}}
\DeclareMathOperator{\agenus}{\mathsf{a-genus}}
\DeclareMathOperator{\ggenus}{\mathsf{g-genus}}
\DeclareMathOperator{\num}{\mathsf{num}}
\DeclareMathOperator{\nnum}{{\mathsf{num}}^\prime}
\DeclareMathOperator{\nnumg}{{\mathsf{num}}^\prime_{\it g}}
\DeclareMathOperator{\nnumm}{{\mathsf{num}}^\prime_{\it m}}
\DeclareMathOperator{\nnumn}{{\mathsf{num}}^\prime_{\it n}}
\DeclareMathOperator{\nnummn}{{\mathsf{num}}^\prime_{{\it m}+{\it n}}}
\DeclareMathOperator{\length}{\mathsf{length}}
\DeclareMathOperator{\im}{\mathsf{Im}}

\DeclareMathOperator{\Hom}{\mathsf{Hom}}
\DeclareMathOperator{\Aut}{\mathsf{Aut}}
\DeclareMathOperator{\Mod}{\mathsf{Mod}}

\def\R{{\mathbb R}}
\def\Z{{\mathbb Z}}

\def\H{{\mathcal H}}

\def\F{{\mathbb F}}
\def\Fa{{\mathbb F_1}}
\def\Fb{{\mathbb F_2}}

\def\B{{\mathcal B}}
\def\S{{S}}

\def\surface{{\mathbb S}}

\def\graph{{\Gamma}}
\def\conjclass{\omega}
\def\MR{\textsf{MR}}
\def\homo{\phi}
\def\homoa{\phi}
\def\homob{h}
\def\homoc{h_2}

\def\circle{C}
\def\tight{\tau}
\def\nv{\mathcal N\mathcal V}
\def\ne{\mathcal N\mathcal E}

\def\d{\partial}
\def\bounding{representation}
\def\primitive{indivisible}
\def\G{H}
\def\P{\mathcal P}

\newcommand\tbound[1]{\hat{#1}}

\author{Mladen Bestvina\and Mark Feighn}
\thanks{Both authors gratefully acknowledge support of the NSF} 

\title[Counting maps]{Counting maps from a surface to a graph}
\begin{document}
\maketitle
\tableofcontents
\section{Introduction and results}
Fix a nonabelian free group $\F$ of finite rank and let $\G$ be a
finitely generated (or f.g.\ for short) group with a f.g.\ subgroup
$P$. In his work on the
Tarski problem, Zlil Sela considers the following question. In how many
ways can a given homomorphism $P\to\F$ be extended to $\G$? Of course
without further restrictions the answer is often infinitely many. He
goes on to define a natural equivalence relation on the set of
extensions (described below in our setting) and obtains the remarkable
result:

\begin{thm}[Sela \cite{zs:tarski3}]\label{t:book3}
Suppose that $\G$ is freely indecomposable rel $P$. There is a number
$N=N(\G,P)$ and a finite set $\mathcal F=\{q_i:\G\to L_i\}$ of proper
quotients so that each homomorphism
$P\to\F$ has at most $N$ equivalence classes of extensions to $\G$
with the property that no element of the equivalence class factors
through an element of $\mathcal F$.
\end{thm}

The set $\mathcal F$ is a {\it factor set for $(\G,P)$}. A
homomorphism from $\G$ to $\F$ with the property that no equivalent
element factors through an element of $\mathcal F$ is {\it solid with
respect to $\mathcal F$}. Not much was known about $N(\G,P)$. For
example, Sela asked whether there was a sequence of examples
$(\G_i,P_i)$ with $\lim N(\G_i,P_i)=\infty$. Our main result is that
there is such a sequence. In fact, in our sequence $\G_i$ will be the
fundamental group of an orientable surface of genus $i$ with $P_i$
representing its boundary and we show that $N(\G_i,P_i)\ge 2^i$.

We now describe our results in more detail. Identify $$\G_g=\langle
a_1,b_1,\cdots,a_g,b_g\rangle,$$ a free group of rank $2g$, with the
fundamental group of a surface $\surface_g$ of genus $g$ and one boundary
component and set
$\d_g=[a_1,b_1]\cdots[a_g,b_g]$ so that $\d_g$ is represented by the
boundary of $\surface_g$.  For $x\in\F$, a {\it genus $g$
  representation of $x$} is a homomorphism $\homob\in
\Hom(\G_g,\F)$ such that $\homob(\d_g)=x$. Set $P_g=\langle\d_g\rangle$.

\begin{definition}\label{d:fdt}
Two genus $g$ \bounding s $\homob$ and $\homob'$ of $x$ are related by a {\it
fractional Dehn twist} if one of the following holds:
\begin{itemize}
\item $\G_{g}=A*_C B$ with $C$ cyclic, $\d_g\in A$, and there is
$z\in\F$ centralizing $\homob(C)$ such that
$\homob'=\homob*(i_z\circ\homob)$ (by which we mean
$\homob'|A=\homob|A$ and $\homob'|B=i_z\circ(\homob|B)$ where $i_z$
denotes conjugation by $z$).
\item 
$\G_{g}=A*_\homo=\langle A,t\mid t^{-1}ct=\homo(c), c\in C\rangle$
where $\homo:C\to C'$ is the bonding isomorphism, $C$ is cyclic, $\d_g\in A$, and 
there is $z\in\F$ centralizing $\homob(C)$ such that
$\homob'|A=\homob|A$ and $\homob'(t)=z\homob(t)$.
\end{itemize}
The equivalence relation
``$\sim$'' on \bounding s of $x$ is generated by
$\homob\sim\homob'$ if $\homob$ and $\homob'$ are related by a
  fractional Dehn twist.
\end{definition}

\begin{remark}\label{r:equivalence}
\begin{enumerate}
\item\label{i:curve} It is a result of Stallings that splittings of
$H_g$ as in Definition~\ref{d:fdt} are all induced by some simple
closed curve $\sigma$ in $\surface_g$. The different items correspond
to whether or not $\sigma$ is separating.\footnote{To prove
(\ref{i:curve}), first resolve the given splitting to find a
collection of pairwise disjoint simple closed curves that induces a
splitting of $H_g$ that can be folded to the given splitting. Then,
show that only in trivial situations is a fold possible; see
\cite{bf:bounding}.}
\item\label{i:dehn} If $z\in\homob(C)$ then there is an automorphism
$\tau$ of $\G_g$ such that $\homob'=\homob\circ\tau$. By
  (\ref{i:curve}), $\tau$ is a (classical) Dehn twist.
\item
Using the trivial splitting
$\G_g=\G_g*_{P_g}P_g$, we see that
$i_z\circ \homob\sim \homob$ where $z$ is a root of $x$ in $\F$.
\item
The group of outer automorphisms of $\G_g$ preserving the conjugacy class of
$\partial_g$ may be identified with the modular group
$\Mod(\surface_g)$, see \cite{hz:mcg}. $\Mod(\surface_g)$ is generated
by Dehn twists \cite{md:mcg,md:englishmcg}. It follows that if
$\homo$ is an automorphism of $\G_g$ fixing $\d_g$ then
$\homob\circ\homo\sim\homob$.
\end{enumerate}
\end{remark}

\begin{example}
Let $\F=\langle u,v\rangle$, $\G_1=\langle a_1,b_1\rangle$, and
$x_m=[u^{m},v]$. For $m,n\in \Z$, the homomorphism
$\homob_{m,n}:\G_1\to\F$ given by $a_1\mapsto u^{m}$ and $b_1\mapsto
vu^n$ is a genus 1 \bounding\ for $x_m$. The homomorphisms $\homob_{m,n}$ and
$\homob_{m,n'}$ are related by fractional Dehn twist whereas they
differ by a Dehn twist iff $n\equiv n'\mod{m}$. In particular,
Theorem~\ref{t:book3} is false if the equivalence relation is defined
only using Dehn twists.
\end{example}

\begin{definition}
For $x\in\F$, a genus $g$ representation is {\it $\sim$-injective} if
all equivalent genus $g$ representations are also injective. Define
$$\num_{g} x$$ to be the number of equivalence classes of
$\sim$-injective genus $g$ representations of $x$ in $\Hom(\G_g,\F)$.
Finally, define
$$f_\F(g)=\sup\{\num_{g} x\mid x\in\F \}.$$
\end{definition} 

If $h\in\Hom(\G_g,\F)$ is $\sim$-injective, then
$h$ is solid with respect to any factor set, and in particular
$$f_\F(g)\le N(\G_g,P_g).$$ That $f_\F(g)$ is finite is a consequence of
Theorem~\ref{t:book3}. In Corollary~\ref{c:independent} we show that
$f_\F$ is independent of $\F$.

It is not hard to see that if $x\in\F$ is a ``generic'' element with a
genus 1 representation, then $\num_{g} x=1$. However, it should
also be reasonable to expect that $f_\F(1)>1$ -- take a ``generic''
map from the genus 2 surface to a graph representing $\F$, then the
element $x\in\F$ represented by the image of the waist curve is
written as $[p,q]$ in two inequivalent ways, giving two
representations. It takes a little bit of work to show that these
representations really are inequivalent. This is the content of
Section \ref{genus 1} and reproduces a result of Lyndon and Wicks
\cite{lw:commutators}\footnote{Thanks to Leo Comerford for pointing us
to this article.}.

For higher genera this conceptual argument fails to show
$f_\F(g)>2$. The reason is that we do not know explicitly the
\MR-diagram\footnote{Some comments are meant for those familiar with
Sela's work on the Tarski problems. The theorems and proofs in this
paper do not depend on such a familiarity.} for the group obtained by
gluing say 3 surfaces with boundary along their boundaries. The only
``obvious'' quotients are obtained by identifying two of the surfaces
or killing the common boundary. To find interesting examples one would
have to show that there are other maximal limit group quotients of
this group, see Remark~\ref{r:interesting}.

However, in Section~\ref{higher genus} we will argue that $f_\F(g)\geq
2^g$. For example, to see $f_\F(2)\geq 4$ we form the ``boundary
connected sum'' of genus 1 examples. Each piece bounds in two ways, so
we expect the sum to bound in 4 ways.

In order to deal with fractional Dehn twists
it is convenient to consider more restrictive representations.

\begin{definition}
An injective
representation $\homob\in\Hom(\G_g,\F)$ of $x\in\F$ is {\it admissible} if
$\im\homob$ is a
primitive\footnote{closed under taking roots, {\it root-closed} in
\cite{lw:commutators}} subgroup of $\F$.
\end{definition}

\begin{prop}
Suppose $\homob\in\Hom(H_g,\F)$ is an admissible \bounding\ for $x\in\F$. If
$\homob'\sim\homob$, then $\homob'$ is also admissible and
$\im \homob'=\im\homob$. In particular, $h$ is $\sim$-injective.
\end{prop}

\begin{proof}
Simple closed curves represent \primitive\footnote{not a proper
power}\ elements of $H_g$, and hence (in the presence of
admissibility) fractional Dehn twists are Dehn twists, see
Remark~\ref{r:equivalence}(\ref{i:dehn}). It follows that there is an
automorphism $\tau$ of $\G_g$ such that $\homob'=\homob\circ\tau$ and
$\homob'(\G_g)=\homob\circ\tau(\G_g)=\homob(\G_g)$.
\end{proof}

\begin{definition}
For $x\in\F$ define $\nnumg x$ to be
$$\big|\{\im\homob:\homob\in\Hom(H_g,\F)\mbox{ is an admissible
representation of }x\}\big|$$ and
$$f_\F'(g)=\sup\{\nnumg x\mid x\in\F\}.$$
\end{definition}

We then have
$$N(\G_g,P_g)\geq f_\F(g)\geq f_\F'(g).$$
We will see that $f'_\F$ is also independent of $\F$. Our main results
are:

\begin{thm}
$f'_\F(1)\geq 2$ and $f'_\F(m+n)\geq f'_\F(m)f'_\F(n)$.
\end{thm}

\begin{cor}
$N(\G_g,P_g)\geq f_\F(g)\geq f'_\F(g)\geq 2^g.$
\end{cor}

We go on to consider a class of representations called {\it minimal}.

\begin{definition}\label{d:minimal}
For $x$ in the commutator subgroup $[\F,\F]$ of $\F$ the {\it
algebraic genus of $x$}, denoted $\agenus x$, is the smallest $g\geq
0$ such that there is $h\in\Hom(H_g,\F)$ with $h(\d_g)=x$. A genus $g$
representation of $x$ is {\it minimal} if $g=\agenus x$. We now make
the same definitions as before, but restrict ourselves to minimal
representations. Define $$\num x$$ to be the number of equivalence
classes of minimal $\sim$-injective representations of $x$, $$\nnum
x=\big|\{\im\homob:\homob\mbox{ is a minimal admissible representation
of }x\}\big|, $$ $$\hat f_\F(g)=\sup\{\num x\mid \agenus x=g\}\mbox{, and}$$
$$\hat f'_\F(g)=\sup\{\nnum x\mid \agenus x=g\}.$$
\end{definition}

Again we have $$N(\G_g,P_g)\geq \hat f_\F(g)\geq \hat f_\F'(g).$$ We
will see that $\hat f_\F$ and $\hat f'_\F$ are independent of
$\F$. Also,

\begin{thm}\label{t:minimal main}
$\hat f'_\F(1)\geq 2$ and $\hat f'_\F(m+n)\geq \hat f'_\F(m)f'_\F(n)$.
\end{thm}

\begin{cor}
$N(\G_g,P_g)\geq \hat f_\F(g)\geq \hat f'_\F(g)\geq 2^g.$
\end{cor}

\section{Labeled graphs}\label{s:labelings}
A reference for this section is \cite{js:folding}.  $\F$ is a
non-abelian free group with fixed finite basis $\B$. The cyclic word
obtained by cyclically reducing the $\B$-word $w$ is denoted
$[[w]]$. There is a 1-1 correspondence between cyclically reduced
cyclic $\B$-words and conjugacy classes of elements of $\F$. If
$x\in\F$, then $[[x]]$ denotes its conjugacy class. We will sometimes
blur the distinction between $\B$-words (or cyclic $\B$-words) and the
elements (or conjugacy classes) that they represent.

Let $R_\B$ denote the wedge of $|\B|$ oriented circles with
fundamental group identified with $\F$. $R_\B$ is an example of a {\it
labeled graph}. More generally, a labeled graph is a connected
non-empty finite graph\footnote{1-dimensional $CW$-complex} $\graph$
together with a combinatorial\footnote{cellular taking open edges
homeomorphically to open edges} map $l:\graph\to R_\B$ called a {\it
labeling}. We consider two labelings $l$ and $l'$ to be the same if,
for each edge $e$, the paths $l|e$ and $l'|e$ are homotopic rel
endpoints. Thus, a labeling is equivalent to a choice of
$u(e)\in\B^{\pm 1}:=\B\sqcup\B^{-1}$ for each oriented edge $e$ of
$\graph$ such that $u(e^{-1})=u(e)^{-1}$ where $e^{-1}$ is the edge
opposite to $e$. A labeling also induces labelings of edge paths in
$\graph$.

If $l:\graph\to R_{\B}$ is an immersion and if $\graph$ has no valence
1 vertices then we say that $l$ or $\graph$ is {\it tight}. A {\it
morphism} of labeled graphs $l_1:\graph_1\to R_\B$ and
$l_2:\graph_2\to R_\B$ is a combinatorial map $f:\graph_1\to\graph_2$
that preserves labels, i.e.\ $l_1=l_2\circ f$. An injective
homomorphism $\homo:\F_1\to \F_2$ induces a cellular map $R_{\B_1}\to
R_{\B_2}$ that immerses each edge. A morphism is obtained by
subdividing edges of $R_{\B_1}$. If $l:\graph\to R_{\B_1}$ is a
labeling then $\homo(l):\homo(\graph)\to R_{\B_2}$ is the induced
labeled graph
$$\graph \overset{l}{\to} R_{\B_1}\to R_{\B_2}.$$ 
Similarly, if
$f:\graph_1\to\graph_2$ is a morphism then there is an induced
morphism $\homo(f):\homo(\graph_1)\to\homo(\graph_2)$.

For a labeling $l:\graph\to R_\B$, $\im\pi_1(l)$ is a well-defined
conjugacy class $\H$ of a subgroup of $\F$ and we say that $l$ is a
{\it labeling} for $\H$ or that $l$ {\it represents} $\H$. There is a
1\nobreakdash-1 correspondence between tight labelings of finite
graphs and conjugacy classes of f.g.\ subgroups of $\F$. A labeling
$l:\graph\to R_\B$ of a finite graph can always be {\it folded} until
it is an immersion, see \cite{js:folding}. Valence one vertices can
then be iteratively pruned until it is tight. Let
$\tight(l):\tight(\graph)\to R_\B$ denote the resulting tight
labeling. This tight labeling is unique unless $\graph$ is
contractible in which case $\tight(\graph)$ will consist of a single
vertex.

Based labeled graphs, i.e.\ labeled graphs with a base point, are also
useful. The definitions in Section~\ref{s:labelings} have analogues if
we allow base points. The base point of the $R_\B$ is its unique
vertex. Of course, labelings automatically take base points to base
points. We require that morphisms do the same. A labeling of a based
labeled graph is tight if it is an immersion and the only valence one
vertex, if any, is the base point. A based labeling $l:(\graph,*)\to
(R_\B,*)$ {\it represents} the subgroup $S$ of $\F_\B$ that is
identified with $\im l_*\subset\pi_1(R_\B,*)$. Without the base point $l$ {\it
represents} the conjugacy class in $\F_\B$ of $S$. If $\graph$ is an
oriented circle with base point, then we also say that $l$ represents
the element $x\in\F_\B$ identified with $l_*([\graph])$ where
$[\graph]\in\pi_1(\graph,*)\cong\mathbb Z$ is the generator determined
by the orientation. Without the base point, $l$ represents the
conjugacy class $[[x]]$ of $x$ in $\F_\B$.  There is a 1-to-1
correspondence between tight based labeled graphs and f.g.\ subgroups
of $\F$. As mentioned above, there is a 1-to-1
correspondence between tight labeled graphs and conjugacy classes of
f.g.\ subgroups of $\F_\B$.

\section{Genus 1}\label{genus 1}
Here $\B=\{ u,v\}$ and so $\F$ is a free group of
rank 2. We use the convention that if $w$ is a $\B$-word then $W$
denotes its inverse.

\begin{prop}[Lyndon-Wicks\cite{lw:commutators}]\label{l:gen1}
$f'_\F(1)\geq 2$. Specifically, if $\homob_1$ is the representation given by
$$u\mapsto uvuvv, v\mapsto UUVU$$
and if $\homoc$ is given by
$$u\mapsto vuvv, v\mapsto UUVUV$$ then $\homob_1$ and $i_u\circ\homoc$
are inequivalent admissible representations for
$$uvuvvUUVUVVu=[\homob_1(u),\homob_1(v)]=i_u([\homoc(u),\homoc(v)]).$$
\end{prop}

The proof of Proposition~\ref{l:gen1} will rely on two lemmas.
\begin{lemma}
$\im\homob_1$ and $\im\homoc$ are not conjugate.
\end{lemma}

\begin{figure}
\scalebox{.8}{\includegraphics{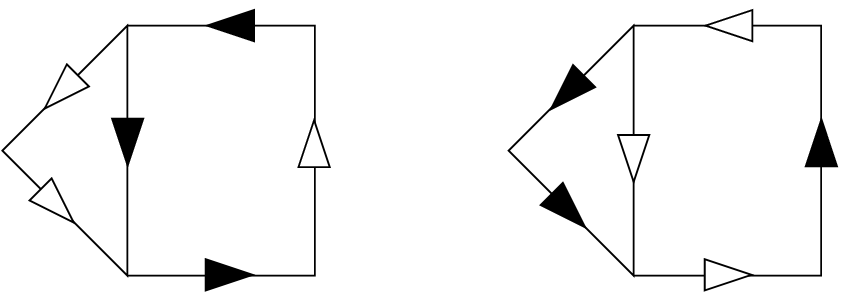}}
\caption{The tight labelings of $\im\homob_1$ and $\im\homoc$.}
\label{f:graphs}
\end{figure}

\begin{proof}
The tight labelings representing the conjugacy classes of $\im\homob_1$
and $\im\homoc$ are pictured in Figure~\ref{f:graphs}. Since these
labelings are not isomorphic, $\im\homob_1$ and $\im\homoc$ are not
conjugate.
\end{proof}

\begin{lemma}
$\im\homob_1$ and $\im\homoc$ are primitive.
\end{lemma}

\begin{proof}
If $\homoa\in\Aut(\F)$ interchanges $u$ and $v$ then
$\homoa(\im\homob_1)=\im\homoc$. So, it is enough to argue that
$\im\homob_1$ is primitive. We will show that $\im\homob_1$ is malnormal in
$\F$, i.e.\ that if $w\in\F$ satisfies 
$i_w(\im\homob_1)\cap\im\homob_1\not=\{1\}$ then $w\in\im\homob_1$. This
clearly implies that $\im\homob_1$ is primitive. The pullback of two
copies of the tight labeling for $\im\homob_1$ has only one component
that is not contractible--that of the ``diagonal''. From
\cite{js:folding}, it follows that $\im\homob_1$ is malnormal in $\F$.
\end{proof}

%\begin{figure}
%\includegraphics{pullback.eps}
%\caption{The off-diagonal components of the pullback.}
%\label{f:pullback}
%\end{figure}

Proposition \ref{l:gen1} is proved.\qed

The homomorphisms $\homob_1$ and $\homoc$ in Proposition~\ref{l:gen1}
were found by a computer search. The original homomorphisms found by
Lyndon and Wicks were $\homob'_1$ given by
$$u\mapsto uvuvUvuvu, v\mapsto vuvuvUvuvUvuvuv$$ 
and $\homoc'$ given by
$$u\mapsto vuvUvuvuvuvuvUvuv, v\mapsto UvuvuvU.$$ 
It is easy to check that $[\homob_1'(u),\homob_1'(v)]$ and
$[\homoc'(u),\homoc'(v)]$ are conjugate. They argue that $\im\homob_1'$
and $\im\homoc'$ are primitive and point out that the abelianizations
of $\homob_1'$ and $\homoc'$ are not in the same $SL_2\Z$-orbit. Hence $\homob_1'$ and
$\homoc'$ are not equivalent.

\section{Higher genus}\label{higher genus}
Here we prove:

\begin{prop}\label{main}
$$f'_\F(m+n)\geq f'_\F(m)f'_\F(n)$$
\end{prop}

\begin{definition}\label{d:alpha map}
Let $\F_1$ and $\F_2$ be two nonabelian free groups with fixed finite
bases $\B_1$ and $\B_2$. For a homomorphism $\homo:\F_1\to\F_2$, set
$m(\homo)=\min\{\length \homo(u)\mid u\in\B_1\}$ where length is with
respect to $\B_2$.\footnote{Recall the convention that an element of
$\F_2$ is identified with the reduced $\B_2$-word representing it.} We
say that $\homo$ is an $\alpha$-map (for some $\alpha>0$) if
\begin{itemize}
%\item for $u,v\in\B_1$ we have
%  $$\frac{\length\homo(u)}{\length\homo(v)}\in
%  (1-\alpha,1+\alpha),$$
\item
for all $u\in\B_1$, a subword of $\homo(u)$ of length
$\ge\alpha m(\homo)$ appears exactly once as a subword of
$\homo(u)$, and
\item for $u,v\in\B_1^{\pm 1}$, if $\homo(u)$ and $\homo(v)$ have
  subwords of length $\ge\alpha m(\homo)$ that are isomorphic preserving
  orientation, then $u=v$.
\end{itemize}
\end{definition}

\begin{remark}
An equivalent definition is that $\homo$ is an {\it $\alpha$-map} if,
for any reduced $\B_2$-word $w$ of length $\ge\alpha m(\homo)$, $w$
appears at most once as a subword in the sequence $\{\homo(u)\mid u\in\B_1^{\pm
1}\}$. For this reason, it is often convenient not to distinguish
between a subword and its inverse. For example, we will say that two
$\B_2$-words $p$ and $q$ {\it share a subword $w$} if $w$ or $W$
appears in $p$ and $w$ or $W$ appears in $q$.\footnote{Recall the convention that corresponding
small and capital letters are mutually inverse.}
\end{remark}

The idea of $\alpha$-maps goes back to Sacerdote \cite{gs:alpha}.

\begin{example}
Say $\F_1=\F_2=\langle u,v\rangle$. Let $$\homo(u)=uvu^2vu^3v\cdots u^nv$$ and
$$\homo(v)=uv^2u^2v^2u^3v^2\cdots u^nv^2$$
As $n\to\infty$, this is an $\alpha$-map for $\alpha\to 0$.
\end{example}

While working with an $\alpha$-map $\homoa:\Fa\to\Fb$ the natural unit
of length is $\alpha m(\homoa)$. We say that an edge path in a
$\B_2$-labeled graph or a $\B_2$-word is {\it $n$-long} if it has
length at least $n\alpha m(\phi)$. Otherwise it is {\it
$n$-short}. %The next principle which follows from the definition of an
%$\alpha$-map will be used repeatedly.

%\vskip 12pt
%\noindent{\it Principle:} Suppose that $u,v\in\B_1^{\pm 1}$, $p$ is a
%subword of $\homoa(u)$, and $p$ shares a 1-long subword with
%$\homoa(v)$. Then $u=v$, $p$ is a subword of
%$\homoa(v)$, and $p$ appears exactly once in
%$\homoa(u)=\homoa(v)$. See Figure~\ref{f:principle}
%%\vskip 12pt

%\begin{figure}
%\includegraphics{principle.eps}
%\caption{The Principle.}
%\label{f:principle}
%\end{figure}

\begin{lemma}\label{l:alpha 1-5} Set $\alpha=1/4$. For all
  $\alpha$-maps $\homo:\F_1\to\F_2$, the following holds.
\begin{enumerate}
\item
$\homo$ is injective. \label{i:alpha injective}
\item\label{i:alpha elements conjugate}
For all $x,x'\in\F_1$, $x$ and $x'$ are $\F_1$-conjugate if and only
if $\homo(x)$ and $\homo(x')$ are
$\F_2$-conjugate.
\item \label{i:alpha element conjugate into subgroup}
For all $x\in\F_1$ and subgroups $\S$ of $\F_1$, $x$ is conjugate into
$\S$ if and only if $\homo(x)$ is conjugate into $\homo(\S)$.
\item\label{i:alpha subgroups conjugate}
For all f.g.\ subgroups $\S$ and $\S'$
of $\F_1$, $\S$ is $\F_1$-conjugate to $\S'$ if and only if
$\homo(\S)$ is $\F_2$-conjugate to $\homo(\S')$.
\item\label{i:alpha primitive element}
For all $x\in\F_1$, $x$ is
\primitive\ in $\F_1$ if and only if $\homo(x)$ is \primitive\ in $\F_2$.
\item\label{i:alpha primitive subgroup} For all subgroups $\S$ of
$\F_1$, $\S$ is primitive in $\F_1$ if and only if $\homo(\S)$ is
primitive in $\F_2$.
\end{enumerate}
\end{lemma}

\begin{proof}
(\ref{i:alpha injective}): Here $\alpha<1/2$ works. Let $l:\circle\to
  R_{\B_1}$ represent the cyclically reduced $\B_1$-word $x=u_1\dots
  u_N$ where $\circle$ is a circle. The labeling $\homo(l)$ represents
  the cyclic $\B_2$-word
  $\homo(x)=\homo(u_1)\cdot{\dots}\cdot\homo(u_N)$ and $\homo(l)$ is nearly
  tight in that folds can only occur in $\alpha
  m(\homo)$-neighborhoods of the initial vertices of the
  $\homo(u_i)$'s. Since $\alpha<1/2$, for each $i$, not all of
  $\homo(u_i)$ is involved in a fold and so $\homo(x)$ is not trivial.

(\ref{i:alpha elements conjugate}): The ``$\Longrightarrow$''
    direction is obvious and holds for any homomorphism
    $\F_1\to\F_2$. For the other direction, assume
    $[[\homo(x)]]=[[\homo(x')]]$. Let $l:\circle\to R_{\B_1}$ be a
    labeling representing $x=u_1\dots u_N$ and let $l':\circle'\to
    R_{\B_1}$ represent $x'=u_1'\dots u'_{N'}$ as cyclically reduced
    cyclic $\B_1$-words. The labelings $\homo(l):\homo(\circle)\to
    R_{\B_2}$ and $\homo(l'):\homo(\circle')\to R_{\B_2}$ represent
    respectively the cyclic $\B_2$-words
    $\homo(u_1)\cdot\homo(u_2)\cdot{\dots}\cdot\homo(u_N)$ and
    $\homo(u'_1)\cdot\homo(u'_2)\cdot{\dots}\cdot\homo(u'_{N'})$. As
    in the proof of (\ref{i:alpha injective}), the labelings
    $\homo(l)$ and $\homo(l')$ are nearly tight. Since $\alpha=1/4$,
    there are 2-long subwords $p_i$ of $\homo(u_i)$ and $p'_j$ of
    $\homo(u_j')$ that survive the folding with $\tight(\homo(l))$ and
    $\tight(\homo(l'))$ representing the same cyclic words $p_1\dots
    p_N=p_1'\dots p'_{N'}$.
\vskip 12 pt
    \noindent{\it Claim:} If $p_i$ and $p_j'$ share a 1-long subword
    $p$ then $p_i=p_j'$.
\vskip 12 pt
\noindent Before proving the claim, we show that it implies
(\ref{i:alpha elements conjugate}). The $p_i$'s and $p'_j$'s are
2-long and so some $p_i$ shares a 1-long subword with some $p_j'$. By
the claim, $p_i=p_j'$. Up to a cyclic permutation, we may assume that
$i=j=1$. Then $p_2$ and $p_2'$ share a 1-long subword and $p_2=p_2'$,
etc.

We now prove the claim. We may assume that $p$ is chosen to be
maximal, i.e.\ $p$ is contained in no longer shared subword. We will
show that $p_i=p=p_j'$. Set $\homo(u_i)=spt$ and
$\homo(u'_j)=s'pt'$. Since $p$ is 1-long, Definition~\ref{d:alpha map}
gives $u_i=u_j'$, $s=s'$, and $t=t'$. Now, $p_i=s_ipt_i$ (so $s_i$ is
the subword of $s$ that survives cancellation). Similarly,
$p_j'=s_j'pt_j'$. The claim is that $s_i$, $s_j'$, $t_i$, and $t_j'$
are all trivial. Since $p$ is maximal one of $s_i$ and $s_j'$, say
$s_i$, is the empty word. If $s_j'$ is not also empty then the
terminal letter of $s'_j$ and the terminal letter of $s$ are the same
letter $b$ and $\homoa(u_{i-1})$ contains the subword $bB$,
contradiction. See Figure~\ref{f:conjugate}. That $t_i$ and $t'_j$
are trivial is similar.

\begin{figure}
\includegraphics{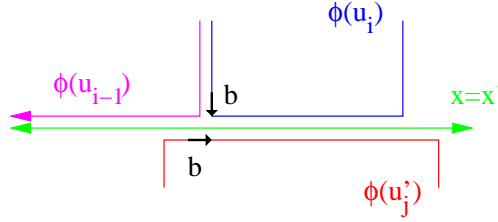}
\caption{
%Lemma~\ref{l:alpha 1-5}(\ref{i:alpha elements
%    conjugate}). 
Adjacent parallel segments should be viewed as overlapping.}
\label{f:conjugate}
\end{figure}

(\ref{i:alpha element conjugate into subgroup}) is a direct
consequence of (\ref{i:alpha elements conjugate}). Indeed, if
$\homo(x)$ is conjugate into $\homo(\S)$ then, for some $s\in \S$,
$\homo(x)$ is conjugate to $\homo(s)$. By (\ref{i:alpha elements
conjugate}) $x$ is conjugate to $s$.

(\ref{i:alpha subgroups conjugate}): Suppose that $\S$ and $\S'$ are
f.g.\ subgroups of $\F_1$ such that $\homo(\S)$ and $\homo(\S')$ are
conjugate in $\Fb$. Let $l:\graph\to R_{\B_1}$ and $l':\graph'\to
R_{\B_1}$ be tight labelings representing the conjugacy classes of
$\S$ and $\S'$ respectively. (\ref{i:alpha elements conjugate}) is a
special case with $\graph=\circle$ and $\graph'=\circle'$. So, we may
assume that $\S$ and $\S'$ are not cyclic.

Consider a natural edge $e$ of $\graph$ viewed as a labeled edge path
representing the word $u_1\dots u_n$. The edge path $\homo(e)$ is a
natural edge of the graph $\homo(\graph)$ representing
$\homo(u_1)\dots\homo(u_N)$. The edge path $\tight(\homo(e))$ nearly
represents a natural edge of $\tight(\homo(\graph))$. That is, there are
2-long subwords $p_i$ of $\homo(u_i)$ so that $p_1\cdots p_N$ is a
natural edge of $\tight(\homo(\graph))$ agreeing with $\tight(\homo(e))$
except perhaps in 1-short initial and terminal subwords. It follows
exactly as in (\ref{i:alpha elements conjugate}) that there is a
corresponding natural edge of $\homo(\graph')$ representing
$\homo(u_1)\cdots\homo(u_N)$ and (\ref{i:alpha subgroups conjugate})
follows.

(\ref{i:alpha primitive element}): The ``$\Longleftarrow$'' direction
is obvious. For the other direction, let $l:\circle\to R_{\B_1}$
represent the non-trivial \primitive\ cyclic word $x=u_1\dots u_N$
which we may assume is cyclically reduced. Suppose that
$\tight(\homo(l)):\tight(\homo(\circle))\to R_{B_2}$ represents
$[[\homo(x)]]=y^n$ with $n>1$ maximal and $y$ cyclically
reduced. Rotation by $2\pi/n$ induces a (label preserving) isomorphism
$\rho:\tight(\homo(\circle))\to\tight(\homo(\circle))$. As in
(\ref{i:alpha elements conjugate}), $y^n=p_1\cdots p_N$ where $p_i$ is
the 2-long subword of $\homo(u_i)$ that survives cancellation. If we
set $p'_i=\rho(p_i)$ then $p_i$ shares a 1-long subword with some
$p'_j$. Exactly as in (\ref{i:alpha elements conjugate}),
$p_i=p'_j$. It follows that $\rho$ leaves the set of $p_i$'s invariant
and that $x$ is not \primitive, contradiction.

(\ref{i:alpha primitive subgroup}) follows directly from (\ref{i:alpha
primitive element}).
\end{proof}

\begin{cor}\label{c:independent}
Set $\alpha=1/4$. If $x\in\Fa$ and if $\homo:\Fa\to\Fb$ is an
$\alpha$-map, then
\begin{itemize}
\item
$\nnumg\homo(x)\geq\nnumg x$; and
\item
$\num_g\homo(x)\geq\num_g x.$
\end{itemize}
In particular, $f'_\F(g)$ and $f_\F(g)$ do not depend on $\F$.
\end{cor}

\begin{proof}
If $\homob$ is an admissible representation of $x$ then
$\homo\circ\homob$ is an admissible representation of $\homo(x)$ by
Items~(\ref{i:alpha injective}) and (\ref{i:alpha primitive subgroup})
of Lemma~\ref{l:alpha 1-5}. Again by Lemma~\ref{l:alpha
1-5}(\ref{i:alpha injective}), $\homo$ induces an injective map from
the set of subgroups of $\F_1$ to the set of subgroups of
$\F_2$. Hence, $\nnum_g\homo(x)\geq\nnum_g x$.

For the second item, let $\homob$ and $\homob'$ be $\sim$-injective
representations of $x$. To show that $\num_g \homo(x)\geq\num_g x$, we must
show two things:
\begin{enumerate}
\item
$\homo\circ\homob$ and $\homo\circ\homob'$ are also $\sim$-injective.
\item
If $\homo\circ\homob\sim\homo\circ\homob'$ then $\homob\sim\homob'$. 
\end{enumerate}

First (1). It is enough to show that $\homo\circ\homob$ is
$\sim$-injective. Suppose $h''$ and $\homo\circ\homob$ are related by
a fractional Dehn twist. There are two cases corresponding to the two
bullets in Definition~\ref{d:fdt}. Suppose $H_g=A*_C B$ as in the
first bullet of Definition~\ref{d:fdt}. Let $z$ be an \primitive\
element of $\F_1$ that commutes with $\homob_i(C)$. By
Lemma~\ref{l:alpha 1-5}(\ref{i:alpha primitive element}), $\homo(z)$
is an \primitive\ element of $\F_2$ that commutes with
$\homo\circ\homob_i(C)$. It follows that for some $k$
$$\homob''=(\homo\circ\homob)*(i_{\homo(z)^k}\circ
(\homo\circ\homob))=\homo\circ (\homob*(i_{z^k}\circ\homob)).$$ The
representation $\homob*(i_{z^k}\circ\homob)$ is injective because it
is equivalent to $\homob$. Since $\homo$ is injective $\homob''$ is
also injective. The case corresponding to the second bullet in
Definition~\ref{d:fdt} is similar and left to the reader.

Continuing with (2), there is a sequence
$$\homob_0,\homob_1,\cdots,\homob_k$$ of representations of $\homo(x)$
where $\homob_0=\homo\circ\homob$, $\homob_k=\homo\circ\homob'$, and
$\homob_{i}$ and $\homob_{i+1}$ are related by a fractional Dehn
twist.  Suppose by induction that $\homob_i=\homo\circ\homob'_i$ for
some $\homob'_i\sim\homob$. Again there are two cases corresponding to
the bullets in Definition~\ref{d:fdt}. Suppose first that the
fractional Dehn twist relating $\homob_i$ and $\homob_{i+1}$ results
from a splitting $H_g=A*_C B$ as in the first bullet of
Definition~\ref{d:fdt}. Let $z$ be an \primitive\ element of $\F_1$
that commutes with $\homob_i(C)$. As in the proof of (1) above,
$\homo(z)$ is an \primitive\ element of $\F_2$ that commutes with
$\homo\circ\homob_i(C)$ and for some $k$
$$\homob_{i+1}=h_i*(i_{\homo(z)^k}\circ h_i)=(\homo\circ
h_i')*(i_{\homo(z)^k}\circ \homo\circ
h'_i)=\homo\circ(h_i'*(i_{z^k}\circ h'_i)).$$ If we set
$h'_{i+1}=h'_i*(i_{z^k}\circ h'_i)$ then
$\homob_{i+1}=\homo\circ\homob'_{i+1}$ and
$\homob'_{i+1}\sim\homob'_i\sim\homob$. At the end of the induction,
we get $\homo\circ\homob'=\homob_k=\homo\circ\homob'_k$. Since $\homo$
is injective, $\homob'=\homob'_k\sim\homob$. Again, the second case is
similar and is left to the reader.

To prove the final statement of this corollary, let $x\in\Fa$ also
satisfy $f_\F(g)=\num_g x$ then $$f_\Fa(g)=\num_g x \le \num_g \homo(x)\le
f_{\Fb}(g).$$ Since $\Fa$ and $\Fb$ were arbitrary,
$f_\Fa(g)=f_\Fb(g)$. The case of $f'_\F$ is similar.
\end{proof}

We are now ready for the proof of our main proposition.

\begin{proof}[Proof of Proposition \ref{main}]
Let $x\in \F$ and $y\in \F$ satisfy $\nnumm x=f'_\F(m)$ and $\nnumn
y=f'_\F(n)$. Since $\nnumm x$ depends only on the conjugacy class of
$x$, we may take $x$ and $y$ to be cyclically reduced. Consider
$z=xy\in\F *\F$. It follows from the next sublemma that $\nnummn z\geq
\nnumm x\cdot\nnumn y$.

\begin{sublemma*}
Suppose that $\homob_x$ and $\homob'_x$ are admissible representations
of $x$ and suppose that $\homob_y$ and $\homob'_y$ are admissible
representations of $y$.
\begin{enumerate}
\item
$\im\homob_x*\homob_y=\im\homob_x*\im\homob_y$ and $\im\homob_x'*\homob_y'=\im\homob_x'*\im\homob_y'$.
\item
$\homob_x*\homob_y:H_{m+n}=H_m*H_n\to\F*\F$ and $\homob'_x*\homob_y'$
are admissible representations of $z$.
\item
If $\im \homob_x\not=\im\homob'_x$ or if
$\im\homob_y\not=\im\homob_y'$ then $\im\homob_x*\homob_y\not=\im\homob'_x*\homob'_y$.
\end{enumerate}
\end{sublemma*}

\begin{proof}
\noindent (1) follows from the uniqueness of normal forms in a
free product; see \cite{ls:book} for example.
\vskip 8pt
\noindent (2): We must show that $\homob_x*\homob_y$ is injective and
has primitive image. Set $A=\im\homob_x$ and $B=\im\homob_y$.
Since the rank of $A*B$ is the sum of the ranks of $A$ and
$B$ and since free groups are Hopfian $\homob_x*\homob_y$ is
injective. Again, \cite{ls:book} is a reference.

By uniqueness of normal forms, an element of $\F*\F$ in normal form
(with respect to $\F*\F$) is in $A*B$ if and only if it is in normal
form with respect to $A*B$. So:
\begin{itemize}
\item an element of $\F*\F$ in
normal form is in $A*B$ if and only if each of its factors
is either in $A$ or $B$.
\end{itemize} Now suppose
$t^n\in A*B$ with $t\in\F*\F$ and $n>0$. We want to show that $t\in
A*B$. To distinguish between the factors of $\F*\F$, let $\F_1$ denote
the first factor and $\F_2$ the second. We want to show $t\in
A*B$. Write $t=t_1t_0^nt_1^{-1}$ as a reduced word in $\F*\F$ with
$t_0$ cyclically reduced. Then, $t^n=t_1t_0^nt_1^{-1}$ is also
reduced. As a first case, suppose that the normal form for $t_0$ has
more than one factor. Let $a$ be the first factor in the normal form
for $t_0$. We may assume that $a\in\F_1$. So, $t_0=aw$ where $w$ is
reduced and has first factor in $\F_2$ and
$t^n=t_1(aw)^nt_1^{-1}$. The occurences of $wa$ in this expression for
$t^n$ are product of factors in the normal form. By the bulleted fact
above, we may remove these factors from $t^n$ and the result is still
in $A*B$. Hence, $t=t_1(aw)t_1^{-1}$ is in $A*B$ and we are done in
this case. The other case is that $t_0\in\F_1\cup\F_2$. We may assume
that $t_0\in\F_1$ and write $t_1=wa$ where $a$ is the last factor of
$t_1$ if it is in $\F_1$ and trivial if the last factor is in
$\F_2$. Using the bulleted fact again, since
$at_0^na^{-1}=(at_0a^{-1})^n$ and $w$ are factors of $t^n$, we
conclude that $at_0^na^{-1}\in A$ and $w\in A*B$. Since $A$ is
primitive, $at_0a^{-1}\in A$ and so $t=wat_0a^{-1}w^{-1}\in A*B$.
\vskip 8pt
\noindent
(3) follows immediately from (1). This finishes the proof of the
sublemma.
\end{proof} 

We now finish the proof of Proposition~\ref{main}. We have established
that $\nnummn z\geq \nnumm x\cdot\nnumn y =f'_\F(m)\cdot
f'_\F(n)$. Also, according to Corollary~\ref{c:independent}, for an
$\alpha$-map $\homo:\F *\F\to\F$ with $\alpha=1/4$, we have
$$\nnummn\homo(z)\geq \nnummn z\geq f'_\F(m)f'_\F(n).$$
Hence, $f'_\F(m+n)\geq \nnummn\homo(z)\geq f'_\F(m)f'_\F(n)$.
\end{proof}

\begin{remark}\label{r:interesting}
We discovered a new limit group quotient that does not factor through
any of the obvious quotients. For example, take $\G$ to be the union of
4 genus 2 surfaces with one boundary component glued along their
boundaries. Take $L$ to be the wedge of two genus two surfaces. Map
$\G\to L$ by sending the common boundary to the product of the two
waist curves, and sending each genus two membrane to the ``boundary
connected sum'' of two halves (there are 4 possible combinations --
use all 4).
\end{remark}

\section{More labeled graphs--boundings}
The rest of the paper is devoted to reproving the main results of this
paper in the context of minimal representations; see
Definition~\ref{d:minimal}.

We now consider the problem of extending a labeling $l:\circle\to
R_\B$ of an oriented circle $\circle$ to a surface. Suppose that we
have a way of pairing up edges of $\circle$ so that paired oriented
edges have the same label in $\B$ and are inconsistently oriented with
respect to the orientation induced by $\circle$. There is a labeling
induced on the quotient graph $\graph$ obtained from $\circle$ by
gluing paired edges and there is an induced morphism
$b:\circle\to\graph$. The morphism $b$ has two key properties:
\begin{enumerate}
\item
$b$ is generically 2-to-1 and generically locally of degree 0, i.e.\
the $b$\,-preimage of an open edge consists of two inconsistently
oriented open edges in $\circle$; and 
\item
the {\it Whitehead graphs} of vertices of $\graph$ are connected.
\end{enumerate}
By (2), we mean the following. The link $Lk_{\graph(b)}(v)$ of a
vertex $v$ of $\graph(b)$ is a union of points, one for each oriented
edge with initial endpoint $v$. For each point $\hat v$ in the
$b$\,-preimage of $v$ there is an induced map $Lk_{C}(\hat v)\to
Lk_{\graph(b)}(v)$.  The {\it Whitehead graph of $v$} has vertex set
$Lk_{\graph(b)}(v)$ and an edge connecting the vertices in the image
of $Lk_{C}(\hat v)\to Lk_{\graph(b)}(v)$ for each $\hat v\in
b^{-1}(v)$. For any $b$ satisfying (1), the Whitehead graph of a
vertex of $\graph$ is a disjoint union of circles. So, to require that
Whitehead graphs are connected is equivalent to requiring them to be
circles. There is a 1-to-1 correspondence between pairings of edges of
$\circle$ as above and morphisms $b$ satisfying (1) and (2). 

\begin{definition}
Suppose that $l:\circle\to R_\B$ is a labeling where $\circle$ is an oriented
circle. A {\it bounding of} $l$ is a morphism $b$ satisfying (1) and
(2) from $l$ to a labeling  $l(b):\graph(b)\to R_\B$.
\end{definition}
\noindent We say that two closed edges of $\circle$ with the same $b$-image are
{\it $b$-paired}. The mapping cylinder $\surface$ of $b$ is a surface
with boundary $C$.  Let $\nv(b)$ denote the set of natural vertices of
$\graph(b)$, i.e.\ the set of vertices of valence other than 2 and let
$\ne(b)$ denote the set of natural edges of $\graph(b)$, i.e.\ the
closures of components of $\graph(b)\setminus\nv(b)$. Set
$v(b)=|\nv(b)|$ and $e(b)=|\ne(b)|$. The {\it geometric genus} of $b$
is defined to be
$$\ggenus\, b=\frac{1}{2}\cdot\big(1-v(b)+e(b)\big)$$ and equals the genus
of $\surface$. If $l$ represents a cyclic $\B$-word $w$ then we also
say that $b$ is a {\it bounding} of $w$ (or of the conjugacy class
$[[w]]$). The {\it geometric genus} of the conjugacy class
$\conjclass$ of an element in $[\F,\F]$ is
$$\ggenus\, \conjclass:=\min \{\ggenus\, b\mid b\mbox{ is a bounding
of } \conjclass\}$$
If $\circle$ is the
concatenation of edge paths $p_1\cdots p_{4g}$ and if the induced edge
paths $b_*(p_j)$ and $b_*(p_{j+2}^{-1})$ coincide for $j\equiv 1$ or
$2\mod 4$, then $b$ is a {\it standard bounding}.

%\hyphenation{re-pre-sent-a-tions}
Of course, there is a close relationship between boundings and
representations. Choose a base point $*\in\circle$ and suppose that
$l$ represents $x\in\F_\B$. Also, choose an isomorphism $$H_g\to
\pi_1(\graph(b), b(*))$$ so that $\d_g$ maps to the generator of $\im
b_*$ determined by the orientation for $\circle$. Here $g=\ggenus
b$. Since we have identified the fundamental group of $R_\B$ with
$\F_\B$, a genus $g$ representation of $x$ is given by:
$$H_g\to\pi_1(\graph(b),b(*))\overset{l(b)_*}{\to}\pi_1(R_\B,*).$$ The
choices here were the base point of $\circle$ and the isomorphism
$H_g\to\pi_1(\graph(b),b(*))$. It follows that if $\homob_1$ and
$\homob_2$ are two representations obtained from $b$ in this manner,
then there is $y\in\F_\B$ and a representation $\homob$ that is
equivalent to $\homob_1$ such that $h_2=i_y\circ h$; see
Remark~\ref{r:equivalence}. We may say that $b$ determines a
representation of the conjugacy class of $x$ that is well-defined up
to equivalence.

\begin{example}\label{e:bounding}
In Figure~\ref{f:bounding}, there are three related boundings. The
first $b$ is a standard bounding of $l:\circle\to R_\B$ where $\B=\{
u,v,w\}$ and $l:\circle\to\R_\B$ represents the cyclic word
$[uv,wU]$. The labeling $\hat b$ is a labeling of $\tight(l)$ and
represents the cyclic word $uvwUVW$. One way to create new boundings
from $b$ is to collapse two edges that are $b$-paired (and then ``pull
apart'' any vertices that may have disconnected Whitehead graph). The
bounding $b'$ is obtained by collapsing the two thicker edges and is a
bounding for $vwUVuW$. Note that $[uv,wU]$ and $uvwUVW$ represent the
same conjugacy class, but $vwUVuW$ represents a different conjugacy
class.
\begin{figure}
\scalebox{.8}{\includegraphics{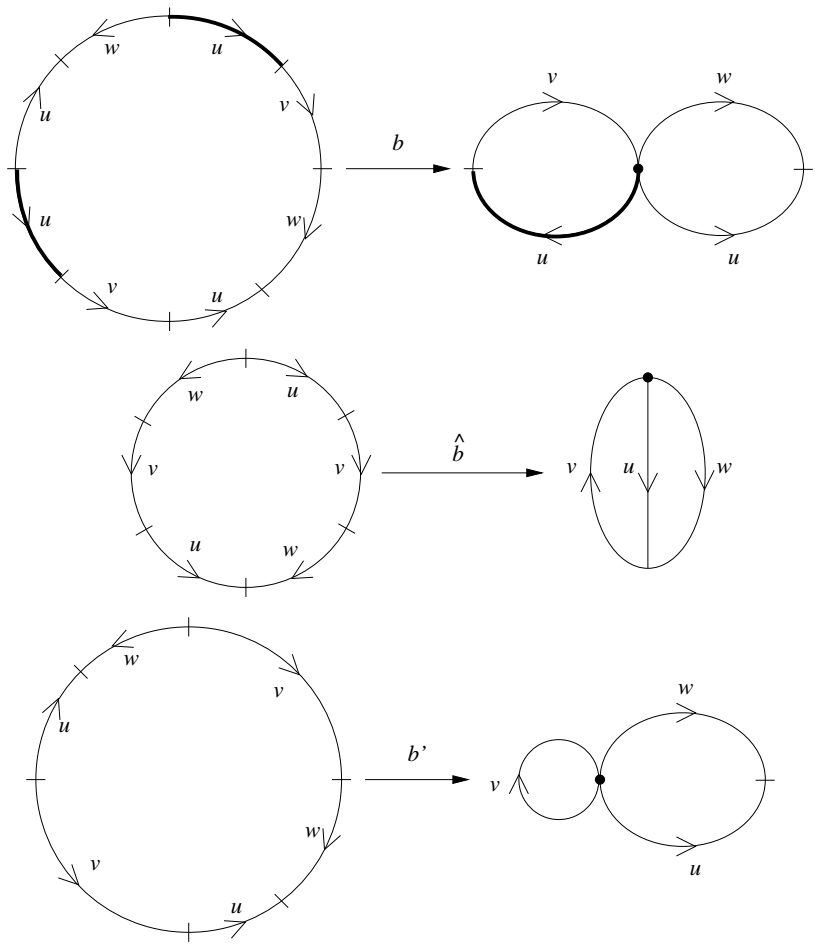}}
\caption{A bounding of a labeled graph, a bounding
  of its tightening, and a bounding obtained by collapsing $b$-paired edges.}
\label{f:bounding}
\end{figure}
\end{example}

\begin{comment}
\begin{example}
The boundings in Figure~\ref{f:2boundings} give two injective
representations of $abaxyxzbabZBABYABAX\in\F_\B$ where
$\B=\{a,b,x,y,z\}$. The first bounding has geometric genus 3 and the
second has geometric genus 2. Although it is not immediately apparent,
the first representation is not $\sim$-injective.\marginpar{This
  example could be safely omitted. Better would be of two
  admissible representations of the same word with different
  genera. Know any?}
\begin{figure}
\scalebox{.7}{\includegraphics{2boundings.eps}}
\caption{Two injective boundings of the same element with different
  geometric genera. The orientation on the circle is clockwise.}
\label{f:2boundings}
\end{figure}
\end{example}
\end{comment}

The next lemma and corollary are classical. The lemma can be proved,
for example, using cut-and-paste surface techniques and folding.

\begin{lemma}\label{l:standard}
Let $b:\circle\to\graph(b)$ be a bounding for the labeling $l:\circle\to
R_{\B}$ representing the cyclic $\B$-word $w$.
\begin{enumerate}
\item Recall that $\tight(l):\tight(\circle)\to
  R_\B$ is the labeling obtained by tightening $l$. There is a
  bounding denoted $\tbound b:\tight(\circle)\to \graph(\tbound b)$ for
  $\tight(l)$ with $\ggenus \tbound b\le\ggenus b$. 
\item
There is a labeled graph $l':\circle'\to R_\B$ representing the
conjugacy class $[[w]]$
with a standard bounding $b':\circle'\to\graph(b')$ such that $\ggenus
b' \le\ggenus b$.\qed
\end{enumerate}
\end{lemma}
\noindent See Figure~\ref{f:bounding}.

\begin{cor}
For $x\in [\F,\F]$, $\agenus x=\ggenus x$.\qed
\end{cor}

\begin{definition}
For $x\in[\F,\F]$, the {\it genus of $x$}, denoted $\genus x$, is the
number $\agenus x=\ggenus x$. Similarly $\genus\, [[x]]:=\agenus\,
[[x]]=\ggenus\, [[x]]$.
\end{definition}

\begin{warning}\label{w:not tight}
The labeled graph $\graph(\tbound b)$ in Lemma~\ref{l:standard}(1)
need not be tight. Even though $\tight(\circle)$ is tight and
therefore $\tbound b$ is an immersion, it is possible that, after a
fold of $\graph(\tbound b)$, the induced map from $\tight(\circle)$ is
no longer generically 2-to-1 and therefore not a bounding. Folding at
bad vertex (see Figure~\ref{f:singularity}) would be an
example. Note however that no folding is possible at a valence two
vertex of $\graph(\tbound b)$.
\end{warning}

We record the next easy lemma for later use.

\begin{lemma}\label{l:easy}
Let $b:\circle\to\graph(b)$ be a bounding for the labeling $l:\circle\to
R_\B$. 
\begin{enumerate}
\item \label{i:collapse} Suppose $b'$ is the new bounding for a new
labeling obtained by first collapsing an edge of $\graph(b)$ and its
$b$\,-preimage and then ``pulling apart'' any vertex with disconnected
Whitehead graph.  Then, $\ggenus\, b'\le\ggenus\, b.$

\item\label{i:bound} 
If $l$ represents a cyclically reduced word then $\graph(b)$ has
no valence one vertices. In particular, $v(b)\le 4\cdot\ggenus\, b-2$
and $e(b)\le 6\cdot\ggenus\, b-3$.\qed
\end{enumerate}
\end{lemma}

The inequalities in (\ref{i:bound}) follow from $2\cdot(\ggenus
b)=1-v(b)+e(b)$ and $3v(b)\le 2e(b)$.

\begin{remark}\label{r:collapse}
The bounding $b'$ in Lemma~\ref{l:easy}(\ref{i:collapse}) is usually a
bounding of a different conjugacy class than the bounding $b$. For
example, see Figure~\ref{f:bounding}.
\end{remark}

\begin{remark}\label{r:confusion}
It is sometimes convenient to view a labeling $l:\graph\to R_\B$ as a
morphism and this can lead to some confusion because the $\homo$-image
of $l$ as a labeling is not usually the same as the $\homo$-image
of $l$ as a morphism. To avoid this confusion, we let $l_\#$ denote the
morphism induced by $l$. See Figure~\ref{f:confusion}.
\begin{figure}
\scalebox{.8}{\includegraphics{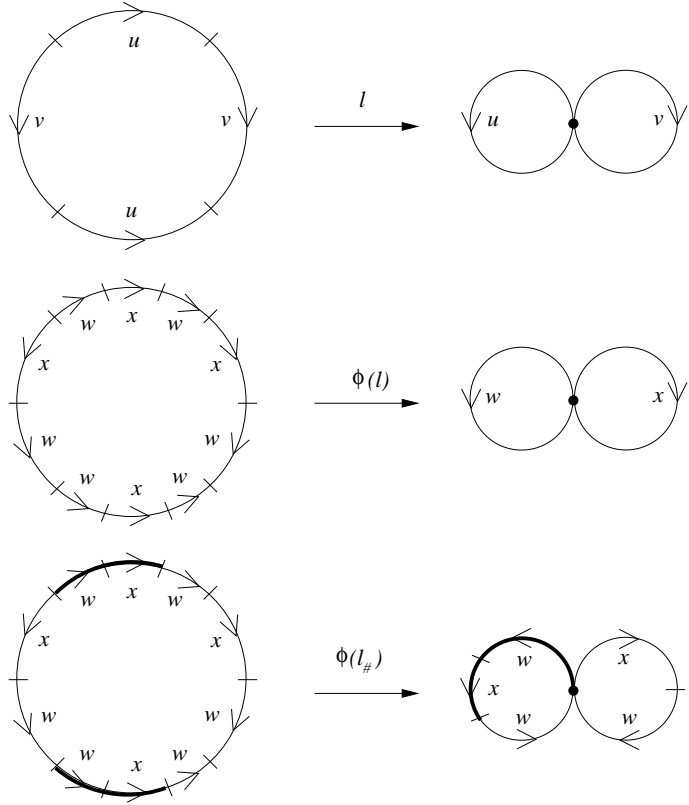}}
\caption{An example of $\homo(l)$ and $\homo(l_\#)$ where
  $\homo:\langle u,v \rangle\to\langle w,x\rangle$ is given by
  $u\mapsto wxw$ and $v\mapsto xw$.}
\label{f:confusion}
\end{figure}
\end{remark}

\section{More $\alpha$-maps--preserving genus}
\begin{lemma}\label{l:alpha 6}
Let $x\in\F_1$ have genus $g$. There is $\alpha>0$ such that, for all
$\alpha$-maps $\homo:\F_{1}\to\F_{2}$, $\homo(x)$ has genus $g$.
\end{lemma}

\begin{proof}
Suppose $x=u_1\cdots u_M\in\F_1$ is cyclically reduced and has genus
$g$. Represent $u_1\cdots u_M$ by a tight labeling $l:\circle\to
R_{\B_1}$ (so $\circle$ has $M$ edges). Choose
$\alpha<[4M(48g-24+M)]^{-1}$. This reason for this choice will become
clear later. Let $\homo:\Fa\to\Fb$ be an $\alpha$-map and set
$m:=m(\homo)$. Consider the induced labeling
$\homo(l):\homo(\circle)\to R_{\B_2}$ (so $\homo(\circle)$ has
$|\homo(u_1)|+\cdots+|\homo(u_M)|$ edges).  We can identify subwords
of $\homo(u_i)$ in $\homo(u_1)\cdots\homo(u_M)$ with certain edge
paths in $\homo(\circle)$. If $w_i$ is a subword of some $\homo(u_i)$
and if $u_i$ equals $u_j$ or $U_j$ then there is a {\it corresponding}
subword $w_j$ of $\homo(u_j)$ or $\homo(U_j)$.  More formally, if
$w_i$ (respectively $w_j$) is represented by the edge path
$p_i:I\to\homo(\circle)$ (respectively $p_j$) then $w_i$ and $w_j$
{\it correspond} if the edge paths $\homo(l_\#)\circ p_i$ and
$\homo(l_\#)\circ p_j$ in $\homo(R_{\B_1})$ are equal, see
Remark~\ref{r:confusion}. The two edge paths in the lower left circle
of Figure~\ref{f:confusion} indicated by the thicker lines correspond.

As in Lemma~\ref{l:alpha 1-5}, $\homo(l)$ is almost tight and
$\tight(\homo(l))$ is obtained by folding $\homo(l)$ in 1-short
neighborhoods of at most $M$ of the vertices of
$\homo(\circle)$. Suppose that $\tight(\homo(l))$ represents the
cyclically reduced word $v_1\cdots v_M$ where each $v_i$ is the
surviving subword of $\homo(u_i)$ (so $\tight(\homo(\circle))$ has
$|v_1|+\cdots+|v_M|$ edges). Since $\alpha<1/4$, the length of each
$v_i$ is at least $m/2$. In order to obtain a contradiction, assume
that $\tight(\homo(l))$ has a bounding $b_{\tight(\homo(l))}$ with
geometric genus $g_{\tight(\homo(l))}$ less than $g$ (see
Lemma~\ref{l:standard}(1)).  Our ultimate goal is to obtain a bounding
for $x$ of geometric genus $\le g_{\tight(\homo(l))}$. By
Lemma~\ref{l:easy}(\ref{i:bound}), $\graph(b_{\tight(\homo(l))})$ has no valence 1
vertices, $v(b_{\tight(\homo(l))})<4g-2$, and
$e(b_{\tight(\homo(l))})<6g-3$. The natural edges of
$\graph(b_{\tight(\homo(l))})$ are labeled with $\B_2$-subwords of
$v_1\cdots v_M$ and, as above, we can talk of their lengths.  We may
also identify the $v_i$'s with edge subpaths of $\homo(\circle)$ via
the labeling $\homo(l)$.  The proof of this lemma will be more
involved than that of Lemma~\ref{l:alpha 1-5} primarily because some
of these natural edges may be 1-short and because
$\graph(b_{\tight(\homo(l))})$ need not be tight (see
Warning~\ref{w:not tight}).  

The edges labeled $u_i$ in $\circle$ map to edge paths labeled
$\homo(u_i)$ in $\homo(\circle)$. Therefore, a bounding of $\homo(l)$
that pairs $\homo(u_i)$'s with $\homo(U_j)$'s can be pulled back to
give a bounding of $l$ with the same geometric genus (forget the
$\homo$'s). Call such a bounding of $\homo(l)$ {\it good}. A weaker
condition on a bounding of $\homo(l)$ is that it be {\it saturated},
i.e.\ paired edges correspond. It is easy to see that a saturated
bounding $b_{\homo(l)}$ of $\homo(l)$ is good if it has the additional
property:
\begin{itemize}
\item
for each natural vertex $v$ of $\graph(b_{\homo(l)})$,
$b^{-1}_{\homo(l)}(v)$ consists only of initial and terminal vertices
of $\homo(u_i)$'s, i.e.\ $b^{-1}_{\homo(l)}(v)$ contains no vertices
that are interior to a $\homo(u_i)$.
\end{itemize} 
Indeed, if this is the case then each natural edge of
$\graph(b_{\homo(l)})$ is a union of $\phi(u_i)$'s and, because
$b_{\homo(l)}$ is saturated $\homo(u_i)$'s will be paired with
$\homo(U_j)$'s.

So, our proof can be completed in two steps. In the first step, we
find a saturated bounding for $\homo(l)$ of geometric genus at most
$g_{\tight(\homo(l))}$. In the second step, we adjust the bounding
discovered in the first step without increasing geometric genus
until it satisfies the bulleted property
above and so is good.

\vskip 12pt
\noindent{\it Step 1.} (Find a saturated bounding $b_{\homo(l)}$ of
${\homo(l)}$ with geometric genus at most $g_{\tight(\homo(l))}$.)
Consider a point $y$ in a natural edge $e$ of
$\graph(b_{\tight(\homo(l))})$ whose distance from
$\nv(b_{\tight(\homo(l))})$ is at least $4\alpha m$. Since the length
of each $v_i$ is more than $m/2$ and $\alpha<1/8$, the
$b_{\tight(\homo(l))}$-image of some $v_j$ meets $e$ in a 2-long
maximal subpath $p$ containing $y$, i.e.\ if we view $v_j$ as a path
in $\graph(b_{\tight(\homo(l)})$ then $p$ is the maximal common
subpath of $v_j$ and $e$ containing $y$. Further, the
$b_{\tight(\homo(l))}$-image of some $V_k$, $k\not= j$ shares a
maximal 1-long subpath $q$ with $p$. Arguing exactly as in
Lemma~\ref{l:alpha 1-5}(\ref{i:alpha elements conjugate}), $p=q$ and
the maximal common subpaths of $v_j$ and $V_k$ (again viewed as paths
in $\graph(b_{\tight(\homo(l))})$) in $e$ and containing $p$
(equivalently $y$) correspond. We conclude that an edge of
$\tight(\homo(\circle))$ whose $b_{\tight(\homo(l))}$-image contains a
point outside the $4\alpha m$-neighborhood of
$\nv(b_{\tight(\homo(l))})$ corresponds with its
$b_{\tight(\homo(l))}$-paired edge.  The number of edges of
$\graph(b_{\tight(\homo(l))})$ in the $4\alpha m$-neighborhood of
$\nv(b_{\tight(\homo(l))})$ is at most $4\alpha m$ times the number of
directions at vertices in $\nv(b_{\tight(\homo(l))})$ which in turn is
at most $4\alpha m\cdot 2e(b_{\tight(\homo(l))})\le 8\alpha m(6g-3)$
by Lemma~\ref{l:easy}(\ref{i:bound}).  Since boundings are generically
2-to-1, the number of edges of $\tight(\homo(\circle))$ not
corresponding with their $b_{\tight(\homo(l))}$-paired edge is at most
$16\alpha m(6g-3)$.

The difference in the number of edges of ${\homo(\circle)}$ and
$\tight({\homo(\circle)})$ is at most $2\alpha m M$. Viewing the edges
of $\tight({\homo(\circle)})$ as edges of ${\homo(\circle)}$ , we have
a paired off corresponding edges of ${\homo(\circle)}$ except for at
most $16\alpha m(6g-3)+2\alpha mM=2\alpha m(48g-24+M)$ edges. So, at
this point we have a {\it partial bounding $\P$} of edges of
$\homo(\circle)$. The bounding is partial in that not all edges of
$\homo(\circle)$ are {\it $\P$-paired} with another edge, such edges
are {\it $\P$-unpaired}. If an edge is $\P$-paired with some other
edge, then we say that $\P$ is {\it defined} on that edge. A partial
bounding that is defined on all edges determines a bounding.

There are at most $2\alpha m(48g-24+M)$ $\P$-unpaired edges in
$\homo(\circle)$ and two edges that are $\P$-paired correspond. From
$\P$ we want to construct a saturated partial bounding where by {\it
saturated} here we mean a partial bounding that in addition to the
property that $\P$-paired edges correspond also has the property that
if an edge is $\P$-unpaired then all corresponding edges are also
$\P$-unpaired. This can be achieved by starting with $\P$ and
forgetting $\P$-pairings of all edges that correspond to a
$\P$-unpaired edge. Since an edge has at most $M$ corresponding edges,
we now have a saturated partial bounding, still called $\P$, of edges
of ${\homo(\circle)}$ that is defined on all but at most $2\alpha
mM(48g-24+M)<m/2$ edges. Since $|v_i|\ge m/2$, in each $v_i$ and
hence in each $\homo(u_i)$, there is at least one edge on which $\P$
is defined. This explains our choice of $\alpha$.

Consider the bounding $b'$ induced from $\P$ by collapsing to a point
each $\P$-unpaired edge of $\homo(\circle)$ as in
Lemma~\ref{l:easy}(\ref{i:collapse}). By construction, two edges that
are $\P$-paired are $b_{\tight(\homo(l))}$-paired. So, $b'$ can also
be obtained by first collapsing to a point each edge of
$\homo(\circle)$ that is not in some $v_i$ (giving
$\tight(\homo(\circle))$) and then iteratively collapsing to points
two $b_{\tight(\homo(l))}$-paired edges that are $\P$-unpaired.  By
Lemma~\ref{l:easy}(\ref{i:collapse}), $\ggenus b'\le
g_{\tight(\homo(l))}$. As noted in Remark~\ref{r:collapse}, $b'$ is
probably not a bounding for $\homo(l)$, but nonetheless we will use
$b'$ and the fact that $\P$ is saturated to complete Step 1 by
extending $\P$ to the sought-after bounding $b_{\homo(l)}$ of
${\homo(l)}$ with $\ggenus b_{\homo(l)}=\ggenus b'\le
g_{\tight(\homo(l))}$.

Recall that $\homo(l)$ represents $\homo(u_1)\dots \homo(u_M)$ and we
may view the $\homo(u_i)$'s as edge paths in $\homo(\circle)$. Suppose
that $p$ is a non-trivial maximal subpath of some $\homo(u_i)$
consisting of $\P$-unpaired edges. Since $\homo(u_i)$ contains an edge
on which $\P$ is defined, an edge $w$ of $p$ shares an endpoint with
an edge $q$ of $\homo(u_i)$ on which $\P$ is defined. Since $\P$ is
saturated, it is defined on all edges of $\homo(\circle)$
corresponding to $q$ and determines a pairing of edges corresponding
to $w$ as follows. If $q_1$ and $q_2$ are $\P$-paired edges
corresponding to $q$ and if $w_k$ corresponds to $w$ and shares an
endpoint with $q_k$, $k=1,2$, then pair $w_1$ with $w_2$. In this way,
we extend $\P$. The extended partial bounding is still saturated and
has fewer unpaired edges. Further, if we now collapse edges that are
unpaired with respect to the extended partial pairing then we get a
bounding $b''$ such that $\graph(b')$ is obtained from $\graph(b'')$
by collapsing disjoint partial natural edges. In particular, $\ggenus
b''=\ggenus b'$. Continue until there are no unpaired edges. This
completes Step~1.

\vskip 12pt
\noindent{\it Step 2.} (Find a good bounding of ${\homo(l)}$ of
geometric genus less than $g$.)  We start with $b_{\homo(l)}$ found in
Step~1. Recall that $b_{\homo(l)}$ is saturated in that
$b_{\homo(l)}$-paired edges correspond.  As previously mentioned, if,
for each natural vertex $v$ of $\graph(b_{\homo(l)})$,
$b_{\homo(l)}^{-1}(v)$ consists only of initial and terminal vertices
of $\homo(u_i)$'s then $b_{\homo(l)}$ would be the desired bounding. A
natural vertex $v$ not having this property is {\it bad}. Let $N(v)$
be the closed neighborhood of $v$ consisting of the union of all
closed edges incident to $v$. We now examine the structure of $N(v)$
for bad $v$.

Suppose $v$ is bad. Give each of the edge paths $\homo(u_i)$ in
$\homo(\circle)$ an orientation so that corresponding $\homo(u_i)$'s
have the same orientation. Since $v$ is bad, each point of
$b^{-1}_{\homo(l)}(v)$ is an interior vertex of some
$\homo(u_i)$. Indeed, because $v$ is bad some element of
$b^{-1}_{\homo(l)}(v)$ is an interior vertex and because
$b_{\homo(l)}$ is saturated all elements are interior vertices. It
follows that $b^{-1}_{\homo(l)}(N(v))$ consists of the vertices
$\tilde v$ in $b^{-1}_{\homo(l)}(v)$ together with, for each $\tilde
v$, the pair of edges $I_{\tilde v}$ and $O_{\tilde v}$ incident to
$\tilde v$. We choose the notation so that, with respect to the
orientation on the $\homo(u_i)$'s, $I_{\tilde v}$ has initial vertex
$\tilde v$ (and so is incoming) and $O_{\tilde v}$ has terminal vertex
$\tilde v$ (and so is outgoing). Finally, since $b_{\homob(l)}$ is
saturated and Whitehead graphs are connected, all $I_{\tilde v}$'s
correspond and all $O_{\tilde v}$'s correspond.

We now introduce a move that produces from $b_{\homo(l)}$ a new
saturated bounding for $\homo(l)$ with no greater geometric genus. We
will then show that we get a good bounding after iterating this move
finitely many times. Intuitively, we ``push the problem
forward''. Since the oriented edges $I_{\tilde v}$ correspond for
$\tilde v\in b^{-1}_{\homo(l)}(v)$, they are all labeled with the same
element $i_v$ of $\B_2$. Similarly all $O_{\tilde v}$'s are labeled
with the same element $o_v$ of $\B_2$. The new bounding is obtained by
collapsing the $O_{\tilde v}$'s, relabeling the $I_{\tilde v}$'s with
$i_vo_v$, and pulling apart any vertices with disconnected Whitehead
graph. See Figure~\ref{f:singularity}. It is clear that the new
bounding has the advertised properties.

For a vertex $\tilde v$ in $\graph(\circle)$ that is interior to some
$\homo(u_i)$, define $|\tilde v|$ to be the distance from $\tilde v$
to the terminal endpoint of $\homo(u_i)$ (remember our orientation on
the $\homo(u_i)$'s). It is easy to check that the following number
decreases upon each iteration: $$\sum \{|\tilde v|:b_{\homo(l)}(\tilde
v)\mbox{ is a bad vertex of }\graph(b_{\homo(l)})\}.$$
This completes Step~2 and the proof of the
lemma.

\begin{figure}
\includegraphics{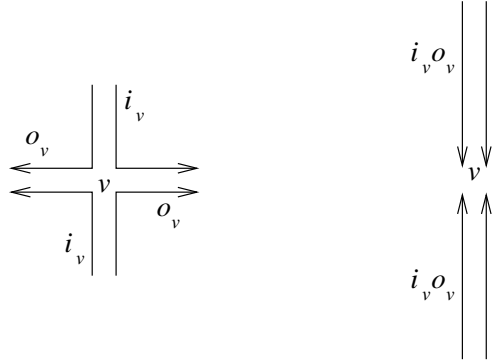}
\caption{On the left are four corresponding paths fitting together to
  form $N(v)$ for a bad vertex $v$ of $\graph(b_{\homo(l)})$. On the
  right is the result of the move described in Step~2 of the proof of
  Lemma~\ref{l:alpha 6}. The
  arrows indicate the orientation induced from the chosen orientation
  of the $\homo(u_i)$'s.}
\label{f:singularity}
\end{figure}
\end{proof}

\section{Proof of main results for minimal representations}
\begin{cor}\label{c:independent2}
Given $x\in\Fa$ there is $\alpha>0$ such that, for any $\alpha$-map
\begin{itemize}
\item
$\nnum\homo(x)\geq\nnum x$, and
\item
$\num\homo(x)\geq\num x.$
\end{itemize}
In particular, $\hat f'_\F(g)$ and $\hat f_\F(g)$ do not depend on $\F$.
\end{cor}

\begin{proof}
The proof is the same as that of
Corollary~\ref{c:independent} as long as we choose $\alpha<1/4$ and
such that $\genus\homo(x)=\genus x$.
\end{proof}

\begin{lemma}\label{l:genus adds}
Suppose that $x$ (respectively $y$) is a cyclically reduced
$\B_1$-word (respectively $\B_2$-word). Let $z=x*y\in\F_{\B_1}*\F_{\B_2}.$
\begin{itemize}
\item
If $b_z:\circle\to R_{\B_1\sqcup\B_2}$ be a bounding for $z$ and let
$b_x$ (respectively $b_y$) be the bounding for $x$ (respectively $y$)
obtained by collapsing to points the edges of $\circle$ labeled with
elements of $\B_2$ (respectively $\B_1$). Then,
$$\ggenus b_x+\ggenus b_y=\ggenus b_z.$$
\item
$\genus x+\genus y=\genus z$
\end{itemize}
\end{lemma}

\begin{proof}
There are two special vertices $c$ and $c'$ in $\circle$ where $x$ and
$y$ meet. Note that $b_z(c)=b_z(c')$. In fact, this is a consequence
of the restriction of $b_z$ to the edge path in $\circle$ labeled
$x$. To see this, complete the edge path labeled $x$ to another circle
$\circle'$ by adding an unlabeled edge (connecting $c$ and $c'$). If
we glue together $b_z$-paired edges of $\circle'$, the quotient is a
surface with boundary and the boundary is the image of the unlabeled
edge. The image of the unlabeled edge is a circle and hence
$b_z(c)=b_z(c')$. It now also follows that Whitehead graphs of
vertices in $\graph(b_x)$ and $\graph(b_y)$ are connected and so $b_x$
and $b_y$ are indeed boundings.

There is a 1-to-1 correspondence between natural vertices in
$\graph(b_z)$ other than $b_z(c)$ and natural vertices of
$\graph(b_x)\sqcup\graph(b_y)$ other than $b_x(c)$ and
$b_y(c)$. Similarly, there is a 1-to-1 correspondence between natural
edges in $\graph(b_z)$ not containing $b_z(c)$ and natural edges of
$\graph(b_x)$ and natural edges of $\graph(b_x)\sqcup\graph(b_y)$ not
containing $b_x(c)$ or $b_y(c)$. Since the four labels of edges
incident to $c$ and $c'$ all have different labels ($x$ and $y$ are
cyclically reduced), $b_z(c)$ is a natural vertex of $\graph(b_z)$. In case
$b_x(c)$ is a natural vertex of $\graph(b_x)$ and $b_y(c)$ is a
natural vertex of $\graph(b_y)$ then there is a 1-to-1 correspondence
between natural edges of $\graph(b_z)$ and natural edges of
$\graph(b_x)\sqcup\graph(b_y)$. In this case,
$$e(b_z)-v(e_z)=e(b_x)-v(e_x)+e(b_y)-v(e_y)-1$$ where the 1 arises
because $b_z(c)$ corresponds to the two vertices $b_x(c)$ and
$b_y(c)$. Hence $\ggenus b_z=\ggenus b_x+\ggenus b_y$. Since the labels in $x$ of the edge incident to and the edge
incident to $c'$ are different, $b_x(c)$ has valence at least two in
$\graph(b_x)$. The same holds for $b_y(c)$. Including vertices of
valence two in our definition of geometric genus produces the same
number. We conclude that in any case $\ggenus b_z=\ggenus b_x+\ggenus
b_y.$ 

The second item follows from the first item and the
observation that if $b_x$ is a bounding for $x$ and if $b_y$ is a
bounding for $y$, then a bounding $b_z$ of $z$ is induced
by ``concatenating'' $b_x$ and $b_y$.
\end{proof}

We can now prove Theorem~\ref{t:minimal main}.

\begin{proof}[Proof of Theorem~\ref{t:minimal main}.]
The representations found in Proposition~\ref{l:gen1} have minimal
genus and hence $\hat f'_\F(1)\ge 2.$ By the second item of
Lemma~\ref{l:genus adds}, the representations found in the proof
of Proposition~\ref{main} are also minimal and so $\hat f'_\F(m+n)\ge
\hat f'_\F(m)\cdot \hat f'_\F(n).$
\end{proof}

\bibliographystyle{plain}
%\bibliography{ref}

\end{document}